\documentclass[12pt,authoryear,preprint]{revtex4}   % style for Physical Review B and AJP are similar

\usepackage{amsmath}    % need for subequations
\usepackage{amsfonts}  %note how statements can be commented out
\usepackage{amssymb}

\begin{document}

\title{Qualitative Criterion for Interception in a Pursuit/Evasion Game}

\author{J. A. MORGAN}
\affiliation{The Aerospace Corporation, P. O. Box 92957 \\Los Angeles, CA 90009, 
United States of America \\
\email{john.a.morgan@aero.org}
}

%\pub{Received (Day Month Year)}{Revised (Day Month Year)}

\begin{abstract}

A qualitative account is given of a differential pursuit/evasion game.  A criterion for 
the existence
of an intercept solution is obtained using future cones that contain all attainable 
trajectories of target or interceptor originating from an initial position.  
A necessary and sufficient conditon that an opportunity to intercept always exist is that, after some 
initial time, the future cone of the target be contained within the future cone of
the interceptor.  The sufficient condition may be regarded as a kind of Nash equillibrium.  

%\keywords{differential games; pursuit/evasion games}
\end{abstract}

%\ccode{MCS Nos.: 49N75, 91A23}  

\maketitle

\section{Introduction}

In this note we consider the differential game that describes the pursuit of a target with position $y(t)$ at time $t$, by an interceptor whose position at time $t$ is $x(t)$.  Both target and interceptor are assumed to maneuver freely and autonomously, subject to certain overall constraints.  The evolution of $x$ and $y$ is given by 
\begin{equation} \frac{dx}{dt}=f(x,u) \label{eq:xeqn} 
\end{equation} 
and 
\begin{equation} \frac{dy}{dt}=g(y,v). \label{eq:yeqn} 
\end{equation} 
Here $f$ and $g$ are assumed to be bounded analytic functions, and $u=u(t)$ and $v=v(t)$ are piecewise analytic controls. The histories of $x$, $y$, and $v$ are assumed to be known up to time $t$. If a control $u$ can be found such that at a time $t_{i}$ the interceptor has contrived to maneuver so that $\parallel x-y \parallel = 0$, an interception is deemed to have taken place, and the game terminates. This problem is a simplified version of one posed by \citet{P1964}.

This study addresses the conditions under which an interception can take place, rather than with conditions of optimality as to time to interception, or payoff functionals.  We seek conditions under which an interception is $\emph{guaranteed}$ to be possible.  There is no assumption, for example, that an interception, if possible, will actually take place.  We are more concerned with the existence of a solution to the game than with payoffs. However, in the present context we may assume that only terminal payoffs are of interest, and that these need not sum to zero. In the terminology of \citet{I1965}, we regard this problem as a game of kind. 

The principal tool used to study this problem is the relation between future cones, as defined in Section (\ref{futurecone}), of target and interceptor.  Constraints on the motion of either player enter the problem formulation through the topology of the future cones. 

One may consider the results proved below as a mathematical exercise confirming the commonplace 
observation in naval warfare that, in order to defeat an adversary in a surface engagement, it is
sufficient and necessary to turn inside the adverary's turning radius. 
Although the discussion will be couched
in terms of trajectories in $\mathbf{R}^{3}$, appropriate for aircraft, missiles, spacecraft, or
simple predator/prey interactions,
%, or critters, 
the results clearly hold in $\mathbf{R}^{2}$, as, for example, in the case of surface naval vessels, 
or in dimensions higher than three.  The extension to more than two players is discussed in the 
sequel.     

\section{The Future Cone} \label{futurecone} 

We begin by defining the future cone of a maneuvering player, and sketching some of its properties. Under the influence of propulsive forces, of its control $u(t)$, and of external forces such as gravity or aerodynamic drag, the position $x(t)$ of the interceptor will evolve as it maneuvers.  The evolution as given by \eqref{eq:xeqn} will be continuous, but need not always be differentiable, depending on the nature of $f$ and of the control $u(t)$. Identical remarks hold for $y(t)$. For an effective interceptor, the control for $x$ will depend upon the history of $y$, in general, and the motion of both interceptor and target may be quite involved.  One may think of their trajectories as being representatives of generic Feynman paths in $\mathbf{R}^{3}$.  But, as a practical matter, both target and interceptor will respect certain physical limitations as to, say, maximum acceleration or total velocity change $\Delta v$.  In addition, the interaction between target and interceptor may be expected to take place during a finite engagement time interval, and to be confined within a finite engagement volume.  

At some initial time $t_{0}$ the interceptor will occupy a definite position $x(t_{0})$.  If we consider all possible histories for $x(t)$ originating at $x(t_{0})$, these will comprise a topological cone in $\mathbf{R}^{4}$ with vertex $x(t_{0})$. We call this set the future cone of $x(t_{0})$ and denote the set of points subsequently accessible to the interceptor in the time interval $(t_{1},t_{2})$, with $t_{0} \le t_{1} < t_{2}$, by$K_{x}^{+}(t_{1}:t_{2};x(t_{0}))$. Clearly, the future cone originating from a vertex at time $t_{1} \subset$ the cone originating from any vertex at a time $t_{0} < t_{1}$.  
%For convenience, the notation $K^{+}_{x}$ or $K^{+}_{x}(t)$ will sometimes be %used when the meaning should be clear from context.

We may suppose that a freely maneuvering player is able to move in any direction in $\mathbf{R}^{3}$ and that, if it can maneuver to a point $x(t)$, it can likewise maneuver to points within any small neighborhood of $x(t)$. It will thus be convenient to regard
the sets $K_{x}^{+}(t_{1}:t_{2};x(t_{0}))$ and $K_{y}^{+}(t_{1}:t_{2};y(t_{0}))$ as open 
for $t_{0} < t_{1} < t_{2}$. (The extension to closed future cones will be discussed in 
Section (\ref{discussion})).The assumption that both target and interceptor can maneuver freely permits us to treat the future cone for each as a path-connected manifold for $t > t_{0}$; in particular, as path-connected for each value of $t$. The attainable trajectories intersecting each point in a spacelike section $K^{+}_{x}(t)$ of a future cone, \emph{i.e.} the subset of the cone at any single value of $t \in (t_{1},t_{2})$, are nowhere tangent to the section, by virtue of the boundedness of $F$.  For $t_{1} \not = t_{2}, K^{+}_{x}(t_{1}) \cap K^{+}_{x}(t_{2}) = \emptyset$. We conclude that the future cones of target and interceptor  for $t > t_{0}$ both admit a timelike foliation \citep{N1967}, and that their spacelike sections are leaves of the foliation.  Every attainable trajectory thus traverses leaves of the future cone in a positive sense as time increases. Given leaf $K^{+}_{x}(t_{0})$, a leaf $K^{+}_{x}(t_{1})$ with $t_{1} > t_{0}$ is generated by exponentiating the action of the tangent space over each point $x(t_{0}) \in K^{+}_{x}(t_{0})$ corresponding to admissible values of $\Delta \mathbf{v}$. The ability of either target or interceptor to maneuver freely also implies that the tangent space over any point in $K^{+}(t_{0})$ is balanced\footnote{The tangent space $T$ over a point is balanced 
if $x \in T$ implies $t \, x \in T$, where $|t| \le 1$}, and thus any leaf to its future will be the union of convex neighborhoods.  As bounded subsets of $\mathbf{R}^{4}$, the sets 
$K_{x}^{+}(t_{1}:t_{2};x(t_{0}))$ and $K_{y}^{+}(t_{1}:t_{2};y(t_{0}))$ have compact closure.
Let $K_{y}^{+}(t_{1}:t_{2};y(t_{0})) \subseteq K_{x}^{+}(t_{1}:t_{2};x(t_{0}))$ and
consider an open covering of $K_{y}^{+}(t_{1}:t_{2};y(t_{0}))$ that includes
$K_{x}^{+}(t_{1}:t_{2};x(t_{0}))$.  Compactness of $\overline{K_{x}^{+}}(t_{1}:t_{2};x(t_{0}))$ thus
implies that $K_{y}^{+}(t_{1}:t_{2};y(t_{0}))$ is covered by a finite union of convex neighborhoods.

The foregoing remarks serve to justify the assumptions we shall make regarding the future 
cones of target and interceptor. It will be assumed that the subset of either cone lying to the future 
of its vertex is a manifold 
%with compact closure 
possessing a timelike foliation, whose leaves are 
locally convex, and which posesses a finite convex cover.  In particular, there is
no assumption that the future cone of either target or interceptor, or their leaves,
are necessarily convex as a whole.

\section{Conditions for interception}

\subsection{Guaranteed Interception}

The existence and character of optimal solutions to equations \eqref{eq:xeqn} and \eqref{eq:yeqn} leading to interception has been well-studied since the work recounted in \citet{P1964}; \emph{vide.} \citep{P1971,P1974,P1981}.  Our concern here is with a qualitative description of conditions under which we may be confident that the interceptor can force an interception, optimally or no.  A winning pure strategy of the interceptor is to choose an attainable trajectory that will lead to interception of the target.  (A mixed strategy would choose among multiple trajectories with this property.) This choice amounts to a mapping $F:K^{+}_{x} \rightarrow K^{+}_{x}$ from the set of all trajectories available to the interceptor to its desired \emph{actual} trajectory. For present purposes, a \emph{guaranteed intercept} will be said to exist when an interception solution always exists, no matter how the target maneuvers within its future cone. For the remainder of this note, by "intercept" is to be understood "guaranteed intercept" even when not explicitly so identified. 
%
%The development in this section relies upon a fixed-point theorem for multifunctions due to 
%Tian~\cite{T1991}.  For completeness, Theorem 3 and Corollary 1 of ref.~\cite{T1991} are 
%restated:
%\newline
%\,
%
%\noindent Theorem 3: \emph{Let $X$ be a nonempty subset of a locally convex separated topological
%vector space $E$.  Suppose that $F:X \rightarrow 2^{E}$ is an upper semicontinuous correspondance with
%nonempty closed convex values.  Then the sufficient and necessary condition for the existence of a 
%fixed point $x^{*} \in F(x^{*})$ is that there exist a nonempty compact convex subset $C \subset X$
%such that} 
%\begin{equation}
%F(x) \cap C \ne \emptyset,\, \forall x \in C.
%\end{equation}
%\newline
%\,
%
%\noindent Corollary 1: \emph{Let $X$ be a nonempty subset of a locally convex separated toplogical 
%vector space $E$.   Suppose that $f:X \rightarrow E$ is a continuous function. Then the sufficient and necessary condition for the existence of a 
%fixed point $x^{*} = f(x^{*})$ is that there exist a nonempty compact convex subset $C \subset X$
%such that} 
%\begin{equation}
%f(x) \cap C, \, \forall x \in C.
%\end{equation}    

\subsection{Interception at a Specified Time:  A Condition on Leaves}

We begin by proving a
\newline
\,

\noindent Lemma:  A necessary and sufficient condition for the existence of a guaranteed intercept is that, at the time of intercept $t_{i}$, \begin{equation} K^{+}_{y}(t_{i};y(t_{1})) \subseteq K^{+}_{x}(t_{i};x(t_{0})). \label{eq:leafcond} \end{equation} for $t_{0}, t_{1} < t_{i}$. 

\emph{Proof}: The necessary condition is elementary.  Suppose a guaranteed intercept exists at time $t_{i}$.  Then, every point $y$ in $K^{+}_{y}(t_{i};y(t_{1}))$ must coincide with some point $x$ in $K^{+}_{x}(t_{i};x(t_{0}))$ in order that $\parallel x-y \parallel = 0$ for at least one pair of values of $x$ and $y$. Were \eqref{eq:leafcond} false, there would be some portion of $K^{+}_{y}(t_{i})$ that lay outside the attainable set of interceptor positions at that time.  Thus there would be a subset of $K^{+}_{y}(t_{i};y(t_{1}))$ for which $\parallel x-y \parallel > 0, \forall x \in K^{+}_{x}(t_{i};x(t_{0}))$.
 
The sufficient condition relies upon a fixed-point theorem  for multifunctions due to \citet{T1991}. Assume that \eqref{eq:leafcond} holds. Then $K^{+}_{x}(t_{i};x(t_{0}))$ is a nonempty compact subset of the separated convex space $\mathbf{R}^{3}$, and its intersection with $K^{+}_{y}(t_{i};y(t_{1}))$ is nonempty. Suppose that, of all the possible trajectories $\subset K^{+}_{y}(t_{i};y(t_{1}))$, the actual trajectory of the target is $y^{*}(t)$.  The mapping from $K^{+}_{x}(t_{i};x(t_{0}))$ into $K^{+}_{y}(t_{i};y(t_{1})) \cap K^{+}_{x}(t_{i};x(t_{0}))=K^{+}_{y}(t_{i};y(t_{1}))$ given by $F(x) = \{y^{*}(t_{i})\}$ is closed and convex, $\forall x \in K^{+}_{x}(t_{i};x(t_{0}))$. For a sequence $x_{n}$ tending to any $x(t_{i}) \in  K^{+}_{x}(t_{i};x(t_{0}))$, $F(x_{n}) = \{y^{*}(t_{i})\}=F(x(t_{i}))$.  The mapping $F$ is thus upper semicontinuous. The point
$y^{*}(t_{i}) \in$ $N_{y^{*}(t_{i})}$ of 
$K^{+}_{x}(t_{i};x(t_{0}))$, where $N_{y^{*}(t_{i})}$ is one of the convex sets that comprise the 
finite open cover of $K_{y}^{+}(t_{1}:t_{2};y(t_{0}))$. Therefore, $F(x)) \cap N_{y^{*}(t_{i})} \not = \emptyset, \forall x \in  N_{y^{*}(t_{i})}.$ The requirements of Theorem 3 of \citet{T1991} are thus satisfied. Combining the consequent fixed point $x^{*}(t_{i})$ for $F$ with the tautological fixed point $y^{*}(t_{i}) \in K^{+}_{y}(t_{i};y(t_{1}))$ resulting from the target's ability to maneuver freely, we may write 
\begin{equation}  \left( \begin{array}{c} x^{*}(t_{i}) \\ y^{*}(t_{i}) \end{array} \right) \in \left( \begin{array}{c} \{y^{*}(t_{i})\} \\ K^{+}_{y}(t_{i};y(t_{1})) \end{array} \right). \label{eq:Nash1} 
\end{equation} $\,\square$ 

\subsection{General Interception Condition} \label{sectIIIB}

We next extend the Lemma to the corresponding assertion for the entire cone $K^{+}_{y}(t_{0}:t_{1};y(t_{1}))$. The conditions for existence of a guaranteed intercept are given by the %\newline \,
\newline
\,

\noindent Theorem: Let $t_{\alpha} < t_{0} < t_{1} < t_{\omega}$. Then, a necessary and sufficient condition for the existence of a guaranteed intercept at $t_{i} \in [t_{0},t_{1}]$ is that \begin{equation} K^{+}_{y}(t_{0}:t_{1};y(t_{0})) \subset K^{+}_{x}(t_{\alpha}:t_{\omega};x(t_{\alpha})) \label{eq:conecond} \end{equation} 

\emph{Proof}: The sufficient condition follows the reasoning used in the Lemma, 
applied to the sets $K^{+}_{y}$ and $K^{+}_{x}$ 
in the locally convex separated topological vector space $\mathbf{R}^{4}$.     
The mapping now is into the target trajectory 
$F(x) = \{y^{*}(t), t \in [t_{0},t_{1}] \}$, which (as before) is closed, convex, and upper 
semicontinuous. To show existence, choose a time 
$t \in [t_{0},t_{1}]$.  By previous remarks, a convex neighborhood $N^{j}_{y^{*}(t)}$
belonging to a finite open cover of $K_{y}^{+}(t_{1}:t_{2};y(t_{0}))$ exists for which
$F(x)) \cap \overline{N}^{j}_{y^{*}(t)} = \{y^{*}(t), t \in [t_{j},t_{j+1}] \} \not = \emptyset$
for some $t_{j},t_{j+1} \in [t_{0},t_{1}] \}$ and $\forall x \in N^{j}_{y^{*}(t)}$.
Then, Theorem 3 of \citet{T1991} gives existence of a fixed point for $t_{0} < t_{i} < t_{1}$
such that
\begin{equation}  
\left( \begin{array}{c} x^{*}(t_{i})\\ y^{*}(t_{i}) \end{array} \right) \in 
\left( \begin{array}{c} \{y^{*}(t), \forall t \in [t_{0},t_{1}] \} \\ 
K^{+}_{y}(t_{0}:t_{1};y(t_{0})) \end{array} \right). 
\label{eq:Nash2}
\end{equation}

The necessary condition is obtained by transfinite induction \citep{K1955}.  Let $t_{\alpha} < t_{0}  < t_{\beta} < t_{\gamma} < t_{i} < t_{1} < t_{\omega}$ and take $t_{i}-t$ as an ordinal.  We prove the result for $K^{+}_{y}(t_{\beta}:t_{1};y(t_{0})) \subset K^{+}_{y}(t_{0}:t_{1};y(t_{0}))$ and extend to the full set $K^{+}_{y}(t_{0}:t_{1};y(t_{0}))$ at the end.

We begin by showing the necessary condition holds at late times.  Suppose that a guaranteed intercept exists at time $t_{i}$ for $t_{\gamma}, t_{i}$ within any neighborhood of $t_{1}$. As $t_{\gamma} \rightarrow t_{1}$, \begin{equation} K^{+}_{y}(t_{\gamma}:t_{1};y(t_{0})) \rightarrow K^{+}_{y}(t_{1};y(t_{0})). \end{equation}  By Lemma 1, it follows that \begin{equation} K^{+}_{y}(t_{1};y(t_{0})) \subseteq K^{+}_{x}(t_{1};x(t_{\gamma})) \subset K^{+}_{x}(t_{\alpha}:t_{\omega};x(t_{\alpha})) \end{equation}
  
Next, suppose that at least one guaranteed intercept opportunity exists for time $t_{i}$ between $t_{\gamma}$ and $t_{1}$.  By the inductive hypothesis, \begin{equation} K^{+}_{y}(t_{\gamma}:t_{1};y(t_{0})) \subset K^{+}_{x}(t_{\alpha}:t_{\omega};x(t_{\alpha})) \end{equation} We wish to examine the prospects at an earlier time $t_{\beta}$.  Consider the sets $K^{+}_{y}(t_{\beta}:t_{\gamma};y(t_{0}))$ and $K^{+}_{x}(t_{\beta}:t_{\gamma};x(t_{\alpha}))$:  \begin{equation} K^{+}_{y}(t_{\beta}:t_{1};y(t_{0}))=K^{+}_{y}(t_{\beta}:t_{\gamma};y(t_{0})) \cup K^{+}_{y}(t_{\gamma}:t_{1};y(t_{0})) \end{equation} and similarly for $K^{+}_{x}$.  If a guaranteed intercept is to be possible  $\forall t \in (t_{\beta},t_{\gamma})$, at no time $t$ in $(t_{\beta},t_{\gamma})$ can it be that \begin{equation}  K^{+}_{y}(t;y(t_{0})) \not \subseteq K^{+}_{x}(t;x(t_{\alpha})), \end{equation} by Lemma 1.  Letting $t_{\beta} \rightarrow t_{0}$, we have \eqref{eq:conecond}. $\,\square$

\section{Discussion} \label{discussion}

The sufficient conditions \eqref{eq:Nash1} and \eqref{eq:Nash2} have been posed in the form of a Nash equillibrium \citep{N1950}. The payoff for the interceptor is positive, while that for the target is negative. One may paraphrase the outcome as follows: The target can navigate to any point in its future cone but, no matter how the target moves, so long as the future cone of the target lies within that for the interceptor, the interceptor can always maneuver to the target's position.

Were $K^{+}_{x}(t_{\alpha}:t_{\omega};x(t_{\alpha}))$ a convex set, the sufficient condition for the
Theorem would
follow immediately from the \citet{K1941} fixed point theorem applied to $\overline{K^{+}_{x}}$.  
This assumption is stronger than might be desired, 
%(the result just obtained is stronger),
but we may suppose the interceptor can 
choose to maneuver within a convex future cone during an engagement, if that is to its advantage. 
This observation amounts to regarding the choice of future cone as part of specifying a strategy for 
the interceptor.

One may treat future cones that are closed sets by taking the closure of the 
sets used in Section (\ref{sectIIIB})\footnote{A closed future cone is equivalent to the 
\emph{attainable set} introduced by \citet{V1967} and \citet{R1969} in a different context.}.
%The target future cone $\overline{K^{+}_{y}}$ has a finite cover of closed convex sets
%$\overline{N}_{y}$. 
Provided that the boundaries of $\overline{K^{+}_{x,y}}$ consist of
trajectories with the properities described in Section (\ref{futurecone}), they satisfy a
transversality condition with respect to the leaves of the cones.  The closed cones thus admit a 
timelike foliation.
%The requirement of transversality of the boundary of a foliated cone 
%$\overline{K^{+}_{x,y}}$ to its leaves is satisfied by any trajectory of $x(t)$, by the remarks in 
%Section (\ref{futurecone}). 
Then, if 
\begin{equation}
\overline{K^{+}_{y}} \subset \overline{K^{+}_{x}},
\end{equation}
both necessary and sufficient conditions for the Theorem are satisfied.  

Extension of these results to dimensionalities other than three poses no difficulties.  As in \citet{N1950}, the extension to an $n$-player game is likewise straightforward.  However, the interpretation of the results differs for the notable case of a single interceptor maneuvering to intercept one genuine target amidst a number of indistinguishable decoys.  Recall that the game terminates when an interception occurs; that is, we do not assume that an interceptor can engage multiple targets in succession. In this case, a guaranteed intercept will certainly exist, but it might lead to the interception of a worthless decoy.  Only if the number of interceptors equals or exceeds the number of targets real and bogus, can one say with confidence that the conditions on the future cones of targets and interceptors obtained in this note guarantee the existence of interception opportunities for all actual targets.

\end{document}